\documentclass[a4paper,12pt]{amsart}
\usepackage{amsmath,amsthm,amssymb,latexsym,epic}
\usepackage{graphicx,enumerate}

\newtheorem{theorem}{Theorem}
\newtheorem{lemma}[theorem]{Lemma}
\newtheorem{remark}[theorem]{Remark}
\newtheorem{corollary}[theorem]{Corollary}
\newtheorem{proposition}[theorem]{Proposition}
\newtheorem{example}[theorem]{Example}

\usepackage[all]{xy}

\begin{document}
\title[Koszul duality for stratified algebras]{Koszul 
duality for stratified algebras I. Quasi-hereditary algebras}
\author{Volodymyr Mazorchuk}
\date{}

\maketitle

\begin{abstract}
We give a complete picture of the interaction between Koszul 
and Ringel dualities for quasi-hereditary algebras admitting linear 
tilting (co)resolutions of standard and costandard modules. We show 
that such algebras are Koszul, that the class of these algebras is 
closed with respect to both  dualities and that on this class these 
two dualities commute. All arguments reduce to short computations in 
the bounded derived category of graded modules.
\end{abstract}

\section{Introduction}\label{s1}

Let $A$ be a positively graded quasi-hereditary algebra. Then there exist
two classical duals for $A$: the Ringel dual $R(A)$ (\cite{Ri}), 
which is the endomorphism algebra of the characteristic tilting $A$-module, 
and the Koszul dual $E(A)$ (\cite{ADL2}), which is the
extension algebra of the direct sum of all simple $A$-modules.
The algebra $R(A)$ is always quasi-hereditary, while the algebra
$E(A)$ is quasi-hereditary only under some additional assumptions.
For example, $E(A)$ is quasi-hereditary if both, projective resolutions 
of all standard $A$-modules and injective  coresolutions of all 
costandard $A$-modules, are linear (see \cite{ADL2}).
Such algebras were called {\em standard Koszul} in \cite{ADL2}.

The natural question to ask is whether $R(E(A))\cong E(R(A))$.
This question was addressed in \cite{MO}, where it was shown that
this is the case under some assumptions, which, roughly speaking,
mean that the algebras $A$, $R(A)$, $E(A)$, $E(R(A))$ and
$R(E(A))$ are standard Koszul with respect to the grading,
induced from the grading on $A$. The main disadvantage of this
result was that the condition was not formulated in terms of 
$A$-modules and hence was very difficult to check.

The main motivation for the present paper was to find an
easier condition which would guarantee the isomorphism
$R(E(A))\cong E(R(A))$. For this we further develop the 
approach of \cite{MO}, based on the category of linear 
complexes of tilting $A$-modules. The main point of the paper is 
that we find an easy way to check Koszulity of $A$ and quasi-heredity
of $E(A)$ based on direct computations in the derived category.
This looks much easier than, for example, the subtle analysis of 
the structure of projective resolutions, carried out in \cite{ADL2}.

A part of the condition, used in \cite{MO}, was formulated as follows: 
all standard and costandard $A$-modules have linear tilting 
(co)resolutions. We call such algebras balanced. Using our 
computational approach we show that 
already this is enough to ensure that all algebras in the list
$A$, $R(A)$, $E(A)$, $E(R(A))$ and $R(E(A))$ are standard Koszul
with respect to the induced grading and derive as a 
corollary that Koszul and Ringel dualities on such $A$ commute.
Under our assumptions we reprove main results from \cite{ADL2} and 
strengthen the main result from \cite{MO}. Our main result is the 
following:

\begin{theorem}\label{thm1}
For every balanced quasi-hereditary algebra $A$ we have:
\begin{enumerate}[(i)]
\item\label{thm1-1} The algebra $A$ is Koszul and standard Koszul.
\item\label{thm1-2} The algebras $A$, $R(A)$, $E(A)$, 
$E(R(A))$ and $R(E(A))$ are balanced.
\item\label{thm1-3} Every simple $A$-module is represented by a linear
complex of tilting modules.
\item\label{thm1-4} $R(E(A))\cong E(R(A))$ as graded 
quasi-hereditary algebras. 
\end{enumerate}
\end{theorem}

By \cite{BGS,MOS} we also have equivalences of the corresponding 
bounded  derived categories of graded  modules for the algebras $A$, 
$E(A)$, $R(A)$ and $R(E(A))\cong  E(R(A))$. Another advantage of 
our approach is that  it admits a straightforward generalization 
to stratified algebras, both in the sense of \cite{ADL1} and \cite{CPS}. 
There is, however, a technical complication in this generalization: 
In the case when a stratified algebra is not quasi-hereditary, it has
infinite global dimension and hence the Koszul dual is
infinite-dimensional. Thus to apply our approach one has first to develop 
a sensible tilting theory for infinite-dimensional stratified algebras. 
This is an extensive technical work, which will be carried out in the 
separate paper \cite{Ma5}. In the present paper we avoid these 
technicalities to make our approach clearer.  Another advantage of our 
approach is that it generalizes to infinite-dimensional 
quasi-hereditary algebras of finite homological dimension.

The paper is organized as follows: In Section~\ref{s2} we collect
all necessary preliminaries about graded quasi-hereditary algebras. 
In Section~\ref{s3} we prove our main result. We complete
the paper with some examples in Section~\ref{s4}.
\vspace{0.5cm}

\noindent
{\bf Acknowledgments.} This research is partially supported by the
Swedish Research Council. Most of the results of the paper were
obtained during the visit of the author to Department of Algebra and
Number Theory, E{\"o}tv{\"o}s University, Budapest in September 2008.
The hospitality of E{\"o}tv{\"o}s University is gratefully acknowledged.
The author also thanks Istv{\'a}n {\'A}goston and Erzs{\'e}bet Luk{\'a}cs
for their hospitality and many stimulating discussions.

\section{Graded quasi-hereditary algebras}\label{s2}

By $\mathbb{N}$ we denote the set of all positive integers.
By a module we always mean a {\em graded left} module, and by grading we
always mean {\em $\mathbb{Z}$-grading}. 
Let $\Bbbk$ be an algebraically
closed field and $A$ be a basic, finite-dimensional, positively graded 
and quasi-hereditary $\Bbbk$-algebra. Let $\Lambda=\{1,\dots,n\}$
and $\{e_{\lambda}:\lambda\in\Lambda\}$ be a complete set of pairwise
orthogonal primitive idempotents for $A$ such that the natural order on 
$\Lambda$ is the one which defines the quasi-hereditary structure on $A$.
Then $A=\oplus_{i\geq 0}A_i$, $A_0\cong \Bbbk e_1\oplus \dots\oplus 
\Bbbk e_n$ and $\mathrm{rad}(A)=\oplus_{i> 0}A_i$. 

Let $A\mathrm{-gmod}$ denote the category of all finite-dimensional 
graded $A$-modules. Morphisms in this category are homogeneous morphism
of degree zero between graded $A$-modules. This is an abelian category
with enough projectives and enough injectives. For $i\in\mathbb{Z}$
we denote by $\langle i\rangle$ the autoequivalence of $A\mathrm{-gmod}$,
which shifts the grading as follows: $(M\langle i\rangle)_j=M_{i+j}$,
$j\in \mathbb{Z}$. We adopt the notation $\mathrm{hom}_A$ and
$\mathrm{ext}^i_A$ to denote homomorphisms and extensions in 
$A\mathrm{-gmod}$.

For $\lambda\in \Lambda$  we consider the graded indecomposable 
projective module $P(\lambda)=Ae_{\lambda}$, its graded simple 
quotient $L(\lambda)=P(\lambda)/\mathrm{rad}(A)P(\lambda)$ and the 
graded indecomposable injective envelop $I(\lambda)$ of $L(\lambda)$.
Let $\Delta(\lambda)$ be the standard quotient of $P(\lambda)$
and $\nabla(\lambda)$ be the costandard submodule of $I(\lambda)$.
By \cite[Corollary~5]{MO}, there exists a graded lift $T(\lambda)$
of the indecomposable tilting module corresponding to $\lambda$ such that 
$\Delta(\lambda)$ is a submodule of $T(\lambda)$ and
$\nabla(\lambda)$ is a quotient of $T(\lambda)$.

For every $i\in \mathbb{Z}$ we will say that {\em centroids} of 
the modules $L(\lambda)\langle i\rangle$, 
$\Delta(\lambda)\langle i\rangle$, $\nabla(\lambda)\langle i\rangle$,
$P(\lambda)\langle i\rangle$, $T(\lambda)\langle i\rangle$
and $T(\lambda)\langle i\rangle$ belong to $-i$. Simple, projective, 
injective, standard, costandard and tilting $A$-modules will be called
{\em structural} modules. A complex $\mathcal{X}^{\bullet}$
\begin{displaymath}
(\mathcal{X}^{\bullet},d_{\bullet}):
\xymatrix{
\dots\ar[r]^{d_{i-2}}&\mathcal{X}^{i-1}\ar[r]^{d_{i-1}}&
\mathcal{X}^{i}\ar[r]^{d_{i}}&\mathcal{X}^{i+1}\ar[r]^{d_{i+1}}&\dots
}
\end{displaymath}
of structural $A$-modules is called {\em linear} provided that 
for every $i\in\mathbb{Z}$ centroids of all indecomposable 
direct summands of $\mathcal{X}^{i}$ belong to  $-i$.

The algebra $A$ is called {\em standard Koszul} provided that all
standard modules have linear projective resolutions and all 
costandard modules have linear injective coresolutions (see \cite{ADL2}).
The algebra $A$ is called {\em balanced} provided that all
standard modules have linear tilting coresolutions and all 
costandard modules have linear tilting resolutions (see \cite{MO}, where
a stronger condition was imposed, however, we will show that both
conditions are equivalent). The algebra $A$ is called {\em Koszul}
provided that projective resolutions of simple $A$-modules
are linear (see \cite{Pr,BGS,MOS}). Denote by $E(A)$ the opposite of 
the Yoneda extension algebra of the direct sum of of all simple 
$A$-modules.  If $A$ is Koszul, the algebra $E(A)$ is called the
{\em Koszul dual} of $A$ and we have that $E(A)$ is Koszul as well
and $E(E(A))\cong A$. 

Let $\mathcal{D}^{b}(A)$ denote the bounded derived category of
$A\mathrm{-gmod}$. For $i\in\mathbb{Z}$ we denote by $[i]$ the 
autoequivalence of $\mathcal{D}^{b}(A)$, which shifts the position
of the complex as follows: $\mathcal{X}[i]^{j}=\mathcal{X}^{i+j}$,
$j\in \mathbb{Z}$ and $\mathcal{X}^{\bullet}\in \mathcal{D}^{b}(A)$. 
As usual, we identify $A$-modules with complexes concentrated in
position $0$. If $A$ is Koszul, then the Koszul duality functor 
\begin{displaymath}
\mathrm{K}=\mathcal{R}\mathrm{hom}_A(\oplus_{i\in\mathbb{Z}}
\mathcal{P}\langle i\rangle[-i]^{\bullet},{}_-),
\end{displaymath}
where $\mathcal{P}^{\bullet}$ is the projective resolution of the direct
sum of simple $A$-modules (see \cite{BGS,MOS}), is well-defined and 
gives rise to an equivalence  from $\mathcal{D}^{b}(A)$ to 
$\mathcal{D}^{b}(E(A))$.

Denote by $\mathfrak{LT}$ the full subcategory of 
$\mathcal{D}^{b}(A)$, which consists of all linear complexes of
tilting $A$-modules. The category $\mathfrak{LT}$ is equivalent to 
$E(R(A))\mathrm{-gmod}$ and the simple objects of $\mathfrak{LT}$ 
have the form $T(\lambda)\langle -i\rangle[i]$, $\lambda\in\Lambda$,
$i\in\mathbb{Z}$ (\cite{MO}). 

Let $R(A)$ denote the {\em Ringel dual} of $A$, which is the opposite
of the (graded) endomorphism algebra of the characteristic tilting 
module $T=\oplus_{\lambda\in\Lambda}T(\lambda)$. The algebra $R(A)$
is quasi-hereditary with respect to the opposite order on $\Lambda$.
The first Ringel duality functor 
\begin{displaymath}
\mathrm{F}=\mathcal{R}\mathrm{hom}_{A}(\oplus_{i\in\mathbb{Z}}T
\langle i\rangle,{}_-)
\end{displaymath}
induces an equivalence from $\mathcal{D}^{b}(A)$ to 
$\mathcal{D}^{b}(R(A))$,  which maps tilting modules to projectives,
costandard modules to standard and injective modules to tilting. 
The second Ringel duality functor 
\begin{displaymath}
\mathrm{G}=\mathcal{R}\mathrm{hom}_{A}({}_-,\oplus_{i\in\mathbb{Z}}T
\langle i\rangle)^*,
\end{displaymath}
where $*$ denotes the usual duality, induces an equivalence 
from $\mathcal{D}^{b}(A)$ to  $\mathcal{D}^{b}(R(A))$,  
which maps tilting modules to injectives,
standard modules to costandard and projective modules to tilting.

\section{The main result}\label{s3}

The aim of this section is to prove Theorem~\ref{thm1}. For this
we fix a balanced algebra $A$ throughout. For $\lambda\in\Lambda$
we denote by $\mathcal{S}_{\lambda}^{\bullet}$ and 
$\mathcal{C}_{\lambda}^{\bullet}$ the linear tilting
coresolution of $\Delta(\lambda)$ and resolution of $\nabla(\lambda)$,
respectively. We will need the following easy observation from 
\cite{MO} and include the proof  for the sake of completeness.

\begin{lemma}[\cite{MO}]\label{lem2}
The natural grading on $R(A)$, induced from $A\mathrm{-gmod}$, 
is positive.
\end{lemma}

\begin{proof}
Let $\lambda,\mu\in\Lambda$. Then $T(\lambda)$ has a standard filtration
and $T(\mu)$ has a costandard filtration (\cite{Ri}). As standard modules 
are  left orthogonal to costandard modules (\cite{Ri}),
every morphism from $T(\lambda)$ to $T(\mu)\langle j\rangle$, 
$j\in\mathbb{Z}$, in induced by a morphism from some standard module
from a standard filtration of $T(\lambda)$ to some costandard module
from a costandard filtration of $T(\mu)$. Hence to prove our claim it is
enough to show that every standard module occurring in the standard
filtration of $T(\lambda)$ and different from $\Delta(\lambda)$ has the 
form $\Delta(\nu)\langle j\rangle$ for some $j>0$; and that every 
costandard module occurring in the costandard filtration of $T(\mu)$ 
and different from $\nabla(\mu)$ has the  form 
$\nabla(\nu)\langle j\rangle$ for some $j<0$.

We will prove the result for $T(\lambda)$ and for $T(\mu)$ the proof
is similar. We use induction on $\lambda$. For $\lambda=1$ the claim
is trivial. For $\lambda>1$ we consider the first two terms of the
linear tilting coresolution of $\Delta(\lambda)$:
\begin{displaymath}
0\to \Delta(\lambda)\to T(\lambda)\to X.
\end{displaymath}
By linearity of our resolution, all direct summands of $X$ have the 
form $T(\nu)\langle 1\rangle$ for some $\nu<\lambda$. All modules from 
the standard filtration of $T(\lambda)$, except for $\Delta(\lambda)$, 
occur in a standard filtration of $X$. Hence the necessary claim follows 
from the inductive assumption.
\end{proof}

From Lemma~\ref{lem2} we directly have the following:

\begin{corollary}\label{cor3}
We have $\mathrm{hom}_A(T(\lambda)\langle i\rangle,T(\mu))=0$,
$\lambda,\mu\in\Lambda$, $i\in\mathbb{N}$.
\end{corollary}

Corollary~\ref{cor3} allows us to formulate the following
main technical tool of our analysis. Let $\mathcal{X}^{\bullet}$ 
and $\mathcal{Y}^{\bullet}$ be two bounded complexes of tilting modules. 
We will say that $\mathcal{X}^{\bullet}$ {\em dominates} 
$\mathcal{Y}^{\bullet}$ provided that for every $i\in\mathbb{Z}$ the 
following holds: if the centroid of an indecomposable summand of 
$\mathcal{X}^{i}$ belongs to $j$ and the centroid of an indecomposable 
summand of $\mathcal{Y}^{i}$  belongs to $j'$, then $j<j'$.  

\begin{corollary}\label{cor4}
Let $\mathcal{X}^{\bullet}$ and $\mathcal{Y}^{\bullet}$ be two 
bounded complexes of tilting modules. Assume that 
$\mathcal{X}^{\bullet}$ dominates $\mathcal{Y}^{\bullet}$.
Then  $\mathrm{Hom}_{\mathcal{D}^{b}(A)}(\mathcal{X}^{\bullet},
\mathcal{Y}^{\bullet})=0$.
\end{corollary}

\begin{proof}
Since tilting modules are self-orthogonal, by 
\cite[Chapter~III(2), Lemma~2.1]{Ha}  the necessary homomorphism
space can be computed already in the homotopy category.
Since $\mathcal{X}^{\bullet}$ dominates $\mathcal{Y}^{\bullet}$,
from Corollary~\ref{cor3} we obtain
$\mathrm{Hom}_{A}(\mathcal{X}^{i},\mathcal{Y}^{i})=0$
for all $i$. The claim follows.
\end{proof}

\begin{proposition}\label{prop5}
For every $\lambda\in\Lambda$ the module $L(\lambda)$ is isomorphic
in $\mathcal{D}^{b}(A)$ to a linear complex 
$\mathcal{L}_{\lambda}^{\bullet}$ of tilting modules.
\end{proposition}

\begin{proof}
Consider a minimal projective resolution $\mathcal{P}^{\bullet}$ of
$L(\lambda)$. Since $A$ is positively graded, for every $i\in\mathbb{Z}$
centroids of all indecomposable projective modules in $\mathcal{P}^{i}$
belong to some $j$ such that $j\geq -i$. Each projective has a
standard filtration. Hence all centroids of standard subquotients in
any standard filtration of an indecomposable projective module 
in $\mathcal{P}^{i}$ also belong to some $j$ such that $j\geq -i$.

Resolving each standard subquotient $\Delta(\lambda)\langle j\rangle$
in every $\mathcal{P}^{i}$ using 
$\mathcal{S}_{\lambda}\langle j\rangle[i]^{\bullet}$, 
we obtain a complex $\overline{\mathcal{P}}^{\bullet}$ of tilting 
modules, which is isomorphic to $L(\lambda)$ in $\mathcal{D}^{b}(A)$. 
By construction and the previous paragraph, for each $i$ all centroids 
of indecomposable summands in $\overline{\mathcal{P}}^{i}$ belong to 
some $j$ such that $j\geq -i$.

Similarly, we consider a minimal injective coresolution
$\mathcal{Q}^{\bullet}$ of $L(\lambda)$.  Since $A$ is positively graded, 
for every $i\in\mathbb{Z}$ centroids of all indecomposable 
injective modules in $\mathcal{Q}^{i}$ belong to some $j$ such that 
$j\leq -i$. Resolving each standard subquotient 
$\nabla(\lambda)\langle j\rangle$ in every $\mathcal{Q}^{i}$ 
using $\mathcal{C}_{\lambda}\langle j\rangle[-i]^{\bullet}$,  we obtain 
another  complex, $\overline{\mathcal{Q}}^{\bullet}$, of tilting modules, 
which is isomorphic to $L(\lambda)$ in $\mathcal{D}^{b}(A)$. By 
construction, for each $i$ all centroids of indecomposable summands in 
$\overline{\mathcal{Q}}^{i}$ belong to some $j$ such that $j\leq -i$.

Because of the uniqueness of the minimal tilting complex
$\mathcal{L}_{\lambda}^{\bullet}$, representing $L(\lambda)$ in 
$\mathcal{D}^{b}(A)$, we thus conclude that for all $i\in\mathbb{Z}$
centroids of all indecomposable summands in 
$\mathcal{L}_{\lambda}^{i}$ belong to $-i$. This means that
$\mathcal{L}_{\lambda}^{\bullet}$ is linear and completes the proof.
\end{proof}

\begin{corollary}\label{cor6}
The algebra $A$ is Koszul.
\end{corollary}

\begin{proof}
Assume that $\mathrm{ext}_A^i(L(\lambda),L(\mu)\langle j\rangle)\neq 0$
for some $\lambda,\mu\in\Lambda$ and $j\in\mathbb{Z}$. Then
$j\leq -i$ as $A$ is positively graded. By Proposition~\ref{prop5},
such a nonzero extension corresponds to a non-zero homomorphism
from $\mathcal{L}_{\lambda}^{\bullet}$ to 
$\mathcal{L}_{\mu}\langle j\rangle[i]^{\bullet}$. Since both
$\mathcal{L}_{\lambda}^{\bullet}$ and 
$\mathcal{L}_{\mu}\langle j\rangle[i]^{\bullet}$ are linear,
the complex $\mathcal{L}_{\lambda}^{\bullet}$ dominates
$\mathcal{L}_{\mu}\langle j\rangle[i]^{\bullet}$ for $j<-i$
and the homomorphism space vanish by Corollary~\ref{cor4}.
Therefore $j=-i$ and the claim follows.
\end{proof}

\begin{corollary}\label{cor7}
The algebra $A$ is standard Koszul.
\end{corollary}

\begin{proof}
That the minimal projective resolution of $\Delta(\lambda)$ is linear,
is proved similarly to Corollary~\ref{cor6}. To prove that the
minimal injective coresolution of $\nabla(\mu)$ is linear we
assume that $\mathrm{ext}_A^i(L(\lambda)\langle j\rangle,\nabla(\mu))\neq 0$
for some $\lambda,\mu\in\Lambda$ and $j\in\mathbb{Z}$. Then
$j\geq i$ as $A$ is positively graded. As both $L(\lambda)$ and
$\nabla(\mu)$ are represented in $\mathcal{D}^{b}(A)$ by linear
complexes of tilting modules, one obtains that for $j>i$ the complex
$\mathcal{L}_{\lambda}\langle j\rangle[-i]^{\bullet}$ dominates
$\mathcal{C}_{\mu}^{\bullet}$, and thus the
extension must vanish by Corollary~\ref{cor4}.
Therefore $j=i$ and the claim follows.
\end{proof}

\begin{corollary}\label{cor8}
The algebra $R(A)$ is balanced.
\end{corollary}

\begin{proof}
By Lemma~\ref{lem2}, the algebra $R(A)$ is positively graded
with respect to the grading, induced from $A\mathrm{-gmod}$.
The functor $\mathrm{F}$ maps linear injective coresolutions of
costandard $A$-modules to linear tilting coresolutions of standard
$R(A)$-modules.  The functor $\mathrm{G}$ maps linear 
projective resolutions of standard $A$-modules to linear tilting 
resolutions of costandard $R(A)$-modules. The claim follows.
\end{proof}

\begin{remark}\label{rem705}
{\rm 
A standard Koszul quasi-hereditary algebra $A$ is balanced if and only
if $R(A)$ is positively graded with respect to the grading induced from
$A\mathrm{-gmod}$, see \cite[Theorem~7]{MO}. 
}
\end{remark}

\begin{corollary}\label{cor805}
The algebra $R(A)$ is Koszul.
\end{corollary}

\begin{proof}
This follows from Corollaries~\ref{cor6} and Corollaries~\ref{cor8}.
\end{proof}

\begin{proposition}\label{prop9}
\begin{enumerate}[(i)]
\item\label{prop9-1} 
The objects $\mathcal{S}_\lambda^{\bullet}$, $\lambda\in\Lambda$,
are standard objects in $\mathfrak{LT}$ with respect to the natural
order on $\Lambda$.
\item\label{prop9-2} 
The objects $\mathcal{C}_\lambda^{\bullet}$, $\lambda\in\Lambda$,
are costandard objects in $\mathfrak{LT}$ with respect to the natural
order on $\Lambda$.
\end{enumerate}
\end{proposition}

\begin{proof}
We prove the claim \eqref{prop9-1}, the claim \eqref{prop9-2} is proved
similarly. Let $\lambda,\mu\in\Lambda$ be such that $\lambda>\mu$. Every
first extension $\xi$ from $\mathcal{S}_\lambda^{\bullet}$ to
$T(\mu)\langle -i\rangle[i]$, $i\in\mathbb{Z}$, is a complex and hence
is obtained as the cone of some morphism $\varphi$ from 
$\mathcal{S}[-1]_\lambda^{\bullet}$ to $T(\mu)\langle -i\rangle[i]$.
The homology of the former complex is $\Delta(\lambda)$ and the
homology of the latter is $T(\mu)$, which has a costandard filtration,
where $\nabla(\lambda)$ does not occur (since $\mu<\lambda$). Since
standard modules are left orthogonal to costandard modules, we get
that all homomorphisms and extensions from $\Delta(\lambda)$
to $T(\mu)$ vanish. Therefore $\varphi$ is homotopic to zero,
which splits $\xi$. The claim follows.
\end{proof}

\begin{proposition}\label{prop10}
For all $\lambda,\mu\in\Lambda$ and $i,j\in\mathbb{Z}$ we have
\begin{equation}\label{eq2}
\mathrm{Hom}_{\mathcal{D}^b(\mathfrak{LT})}
(\mathcal{S}_\lambda^{\bullet},\mathcal{C}_\mu
\langle j\rangle[-i]^{\bullet})=
\begin{cases}
\Bbbk,& \lambda=\mu, i=j=0;\\
0,& \text{otherwise}.
\end{cases}
\end{equation}
\end{proposition}

\begin{proof}
Via the equivalence $\mathrm{K}\circ \mathrm{F}$,
the equality \eqref{eq2} reduces to the equality
\begin{displaymath}
\mathrm{Hom}_{\mathcal{D}^b(A)}
(\Delta(\lambda)^{\bullet},\nabla(\mu)
\langle j\rangle[-i]^{\bullet})=
\begin{cases}
\Bbbk,& \lambda=\mu, i=j=0;\\
0,& \text{otherwise}.
\end{cases}
\end{displaymath}
The latter equality is true as standard modules are left orthogonal 
to costandard modules (see \cite{Ri}).
\end{proof}

\begin{corollary}\label{cor11}
The algebra $E(R(A))$ is quasi-hereditary with respect to the natural
order on $\Lambda$.
\end{corollary}

\begin{proof}
By Propositions~\ref{prop9} and \ref{prop10}, standard $E(R(A))$-modules
are left orthogonal to costandard. Now the claim follows from
\cite[Theorem~1]{DR} (or \cite[Theorem~3.1]{ADL1}).
\end{proof}

\begin{corollary}\label{cor12}
The complexes $\mathcal{L}_{\lambda}^{\bullet}$, $\lambda\in\Lambda$, 
are tilting objects in $\mathfrak{LT}$.
\end{corollary}

\begin{proof}
Because of \cite[Theorem~3.1]{ADL1} (or \cite{DR,Ri}), we just need to show 
that any first extension from a standard object to 
$\mathcal{L}_{\lambda}^{\bullet}$ splits, and that any first extension from 
$\mathcal{L}_{\lambda}^{\bullet}$ to a costandard object splits. We prove the 
first claim and the second one is proved similarly.

Any first extension $\xi$ from 
$\mathcal{S}_{\mu}\langle -i\rangle[i]^{\bullet}$, $\mu\in\Lambda$,
$i\in\mathbb{Z}$, to $\mathcal{L}_{\lambda}^{\bullet}$ is a cone of some
homomorphism $\varphi$ from  
$\mathcal{S}_{\mu}\langle -i\rangle[i-1]^{\bullet}$
to $\mathcal{L}_{\lambda}^{\bullet}$. Thus $\varphi$ corresponds to a 
(nonlinear) extension of degree $1-i$ from $\Delta(\mu)\langle -i\rangle$ to 
$L(\lambda)$. As $A$ is standard Koszul by Corollary~\ref{cor7},
we get that  $\varphi$ is homotopic to zero, 
and thus the extension $\xi$ splits. The claim follows.
\end{proof}

\begin{corollary}\label{cor14}
There is an isomorphism $E(A)\cong R(E(R(A)))$ of graded algebras, both considered with respect to the natural grading induced from 
$\mathcal{D}^b(A)$. In particular, we have $R(E(A))\cong E(R(A))$.
\end{corollary}

\begin{proof}
By Corollary~\ref{cor12}, the algebra $R(E(R(A)))$ is the opposite
of the endomorphism
algebra of $\oplus_{\lambda\in\Lambda}\mathcal{L}_{\lambda}^{\bullet}$.
Since $\mathcal{L}_{\lambda}^{\bullet}$ is isomorphic to 
$L(\lambda)$ in $\mathcal{D}^b(A)$, from 
\cite[Chapter~III(2), Lemma~2.1]{Ha} it follows that the same algebra
is isomorphic to $E(A)$. The claim follows.
\end{proof}

\begin{corollary}\label{cor15}
Both $E(A)$ and $R(E(A))$ are positively graded with respect to
the natural grading induced from  $\mathcal{D}^b(A)$. 
\end{corollary}

\begin{proof}
For $E(A)$ the claim is obvious. By Corollary~\ref{cor14}, we have
$R(E(A))\cong E(R(A))$. As $R(A)$ is positively graded with respect to the
grading induces from $\mathcal{D}^b(A)$ (Lemma~\ref{lem2}), the algebra
$E(R(A))$ is positively graded with respect to the induces grading as well.
\end{proof}

\begin{proposition}\label{prop16}
The positively graded algebras $E(A)$ and $R(E(A))$ are balanced.
\end{proposition}

\begin{proof}
Because of Corollary~\ref{cor8}, it is enough to prove the claim for
the algebra $E(A)$. Consider the algebra $E(R(A))$, whose module
category is realized via $\mathfrak{LT}$. 

\begin{lemma}\label{lem17}
The algebra $E(R(A))$ is standard Koszul.
\end{lemma}

\begin{proof}
We already know that $E(R(A))$ is
positively graded with respect to the grading, induced from 
$\mathcal{D}^b(A)$. Let us show that projective resolutions of 
standard  $E(R(A))$-modules are linear. For injective resolutions
of costandard modules the argument is similar.

We have to compute
\begin{equation}\label{eq3}
\mathrm{hom}_{\mathcal{D}^b(\mathfrak{LT})}
(\mathcal{S}_{\lambda}^{\bullet},T(\mu)\langle j\rangle[i])
\end{equation}
for all $\lambda,\mu\in\Lambda$ and $i,j\in\mathbb{Z}$. Via the
equivalence $\mathrm{K}\circ\mathrm{F}$, the space \eqref{eq3}
is isomorphic to the space
$\mathrm{hom}_{\mathcal{D}^b(A)}
(\Delta(\lambda),T(\mu)\langle j\rangle[i])$.
As $T(\mu)$ has a costandard filtration and standard modules are
left orthogonal to costandard, we get that the later space is non-zero
only if $i=0$. As $R(A)$ is positively graded, we also get that  $j<0$.
Applying \cite[Theorem~22]{MOS} we obtain that the standard 
$E(R(A))$-module $\mathcal{S}_{\lambda}^{\bullet}$ has only linear extensions
with simple $E(R(A))$-modules. This completes the proof.
\end{proof}

Using Lemma~\ref{lem17}, the proof of Proposition~\ref{prop16}
is completed similarly to the proof of Corollary~\ref{cor8}.
\end{proof}

\begin{proof}[Proof of Theorem~\ref{thm1}.]
The claim \eqref{thm1-1} follows from Corollaries~\ref{cor6} and \ref{cor7}.
The claim \eqref{thm1-2} follows from Corollary~\ref{cor8} and 
Proposition~\ref{prop16}. 
The claim \eqref{thm1-3} follows from Proposition~\ref{prop5}. 
Finally, the claim \eqref{thm1-4} follows from Corollary~\ref{cor14}.
\end{proof}

\section{Examples}\label{s4}

\begin{example}\label{exm21}
{\rm
Graded quasi-hereditary algebras, associated with blocks of the usual
BGG category $\mathcal{O}$ and the parabolic category $\mathcal{O}$
for a semi-simple complex finite-dimensional Lie algebra, are 
balanced by \cite{Ma}.
}
\end{example}

\begin{example}\label{exm22}
{\rm
The algebra $A$ is called directed if either all standard 
or all costandard $A$-modules are simple (this is equivalent to
the requirement that the quiver of $A$ is directed with respect to
the natural order on $\Lambda$). For a directed algebra $A$
tilting modules are either injective (if standard modules are simple) 
or projective (if costandard modules are simple). Hence 
any directed Koszul algebra is balanced.
}
\end{example}

\begin{example}\label{exm23}
{\rm
Finite truncations $V_{\mathcal{T}}$ of Cubist algebras from 
\cite[Section~6]{CT} are balanced. Indeed, $V_{\mathcal{T}}$ is
standard Koszul by \cite[Proposition~46]{CT}, and that the Ringel dual
of $V_{\mathcal{T}}$ is positively graded with respect to the induced
grading follows from \cite[Corollary~71]{CT}. So, the fact that 
$V_{\mathcal{T}}$ is balanced follows from Remark~\ref{rem705}. 
}
\end{example}

\begin{example}\label{exm24}
{\rm
One explicit example. Consider the path algebra $A$ of the 
following quiver with relations:
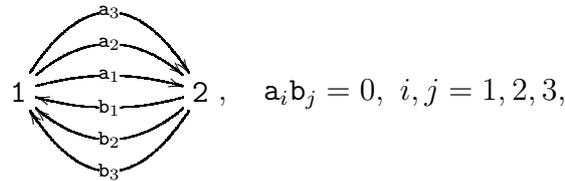
\begin{displaymath}
\xymatrix{
\mathtt{1}
\ar@/^/[rr]|-{\mathtt{a}_1}
\ar@/^1.5pc/[rr]|-{\mathtt{a}_2}
\ar@/^2.5pc/[rr]|-{\mathtt{a}_3}
&&
\mathtt{2}
\ar@/^/[ll]|-{\mathtt{b}_1}
\ar@/^1.5pc/[ll]|-{\mathtt{b}_2}
\ar@/^2.5pc/[ll]|-{\mathtt{b}_3}
},\quad
\mathtt{a}_i\mathtt{b}_j=0,\,\, i,j=1,2,3,
\end{displaymath}
We have $\Delta(\mathtt{1})\cong T(\mathtt{1})\cong
L(\mathtt{1})\cong \nabla(\mathtt{1})$ and for $\lambda=\mathtt{2}$ 
we have the following standard and tilting modules:
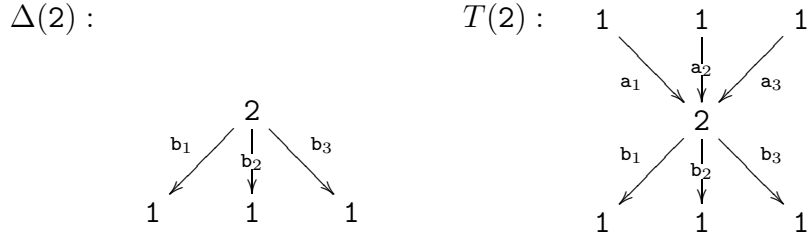
\begin{displaymath}
\Delta(\mathtt{2}):\quad
\xymatrix{\\
&\mathtt{2}\ar[dl]_{\mathtt{b}_1}\ar[d]|-{\mathtt{b}_2}
\ar[dr]^{\mathtt{b}_3}&\\
\mathtt{1}&\mathtt{1}&\mathtt{1}
}\quad\quad\quad
T(\mathtt{2}):\quad
\xymatrix{
\mathtt{1}\ar[dr]_{\mathtt{a}_1}&\mathtt{1}\ar[d]|-{\mathtt{a}_2}
&\mathtt{1}\ar[dl]^{\mathtt{a}_3}\\
&\mathtt{2}\ar[dl]_{\mathtt{b}_1}\ar[d]|-{\mathtt{b}_2}
\ar[dr]^{\mathtt{b}_3}&\\
\mathtt{1}&\mathtt{1}&\mathtt{1}
}
\end{displaymath}
Hence we have the following linear tilting coresolution of 
$\Delta(\mathtt{2})$:
\begin{displaymath}
0\to \Delta(\mathtt{2})\to T(\mathtt{2})\to
T(\mathtt{1})\langle 1\rangle\oplus 
T(\mathtt{1})\langle 1\rangle\oplus
T(\mathtt{1})\langle 1\rangle\to 0.
\end{displaymath}
Swapping $\mathtt{a}_i$ and $\mathtt{b}_i$, $i=1,2,3$, defines an
antiinvolution on $A$, which preserves the primitive idempotents. Hence
there is a duality on $A\mathrm{-gmod}$, which preserves isomorphism
classes of simple modules. Applying this duality to the above resolution
gives a linear tilting resolution of $\nabla(\mathtt{2})$. Thus
$A$ is balanced. In this example one can also arbitrarily increase 
or decrease the number of arrows.
}
\end{example}

\begin{example}\label{exm25}
{\rm
One computes that the path algebra of the 
following quiver with relations
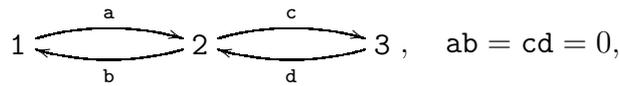
\begin{displaymath}
\xymatrix{
\mathtt{1}
\ar@/^/[rr]^{\mathtt{a}}
&&
\mathtt{2}
\ar@/^/[ll]^{\mathtt{b}}\ar@/^/[rr]^{\mathtt{c}}
&&
\mathtt{3}\ar@/^/[ll]^{\mathtt{d}}
},\quad
\mathtt{a}\mathtt{b}=\mathtt{c}\mathtt{d}=0,
\end{displaymath}
is standard Koszul but not balanced. In fact, the Ringel dual of
this algebra is the path algebra of the following
quiver with relations
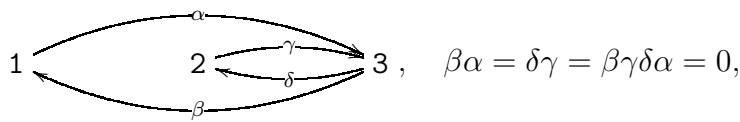
\begin{displaymath}
\xymatrix{
\mathtt{1}
\ar@/^1.5pc/[rrrr]|-{\alpha}
&&
\mathtt{2}\ar@/^/[rr]|-{\gamma}
&&
\mathtt{3}\ar@/^1.5pc/[llll]|-{\beta}\ar@/^/[ll]|-{\delta}
},\quad
\beta\alpha=\delta\gamma=\beta\gamma\delta\alpha=0,
\end{displaymath}
which is not Koszul (not even quadratic). So, our results can not
be extended to all standard Koszul algebras.
}
\end{example}

\begin{remark}\label{rem31}
{\rm
Directly from the definition it follows that if the algebra
$A$ is balanced, then the algebra $A/Ae_n A$ is balanced as well.
It is also easy to see that if $A$ and $B$ are balanced, then
both $A\oplus B$ and $A\otimes_{\Bbbk}B$ are balanced.
}
\end{remark}

\vspace{1cm}

\noindent
Department of Mathematics, Uppsala University, SE 471 06,
Uppsala, SWEDEN, e-mail: {\tt mazor\symbol{64}math.uu.se}

\end{document}